\documentclass[12pt,leqno]{amsart}
\usepackage{amssymb,ifthen,times}
\usepackage[utf8]{inputenc}
\usepackage[T1]{fontenc}
\def\N{\mathbb{N}}
\def\R{\mathbb{R}}
\def\Q{\mathbb{Q}}

\newtheorem{theorem}{Theorem}
\newtheorem*{theorem*}{Theorem}
\def\Thm#1#2{\ifthenelse{\equal{#1}{*}}{\begin{theorem*}#2\end{theorem*}}
             {\begin{theorem}\label{T#1}#2\end{theorem}}}
\newtheorem{Atheorem}{Theorem}

\newtheorem{proposition}[theorem]{Proposition}
\newtheorem*{proposition*}{Proposition}
\def\Prp#1#2{\ifthenelse{\equal{#1}{*}}{\begin{proposition*}#2\end{proposition*}}
             {\begin{proposition}\label{P#1}#2\end{proposition}}}

\newtheorem{corollary}[theorem]{Corollary}
\newtheorem*{corollary*}{Corollary}
\def\Cor#1#2{\ifthenelse{\equal{#1}{*}}{\begin{corollary*}#2\end{corollary*}}
             {\begin{corollary}\label{C#1}#2\end{corollary}}}

\newtheorem{lemma}[theorem]{Lemma}
\newtheorem*{lemma*}{Lemma}
\def\Lem#1#2{\ifthenelse{\equal{#1}{*}}{\begin{lemma*}#2\end{lemma*}}
             {\begin{lemma}\label{L#1}#2\end{lemma}}}
\def\lem#1{Lemma~\ref{L#1}}

\theoremstyle{definition}
\newtheorem{remark}[theorem]{Remark}
\newtheorem*{remark*}{Remark}
\def\Rem#1#2{\ifthenelse{\equal{#1}{*}}{\begin{remark*}\rm #2\end{remark*}}
             {\begin{remark}\label{R#1}\rm #2\end{remark}}}

\newtheorem{example}[theorem]{Example}
\newtheorem*{example*}{Example}
\def\Exa#1#2{\ifthenelse{\equal{#1}{*}}{\begin{example*}\rm #2\end{example*}}
             {\begin{example}\label{Ex#1}\rm #2\end{example}}}

\def\eq#1{{\rm(\ref{E#1})}}
\def\Eq#1#2{\ifthenelse{\equal{#1}{*}}
  {\begin{equation*}\begin{aligned}#2\end{aligned}\end{equation*}}
  {\begin{equation}\begin{aligned}\label{E#1}#2\end{aligned}\end{equation}}}

\begin{document}

\date{\today}

\title[]{An elementary proof for the decomposition theorem of Wright convex functions}

\author[Zs. P\'ales]{Zsolt P\'ales}
\address{Institute of Mathematics, University of Debrecen, 
H-4032 Debrecen, Egyetem t\'er 1, Hungary}
\email{pales@science.unideb.hu}

\subjclass[2000]{Primary 26D15}
\keywords{Wright convexity; Jensen convexity; decomposition theorem.}

\thanks{The research was supported by the EFOP-3.6.1-16-2016-00022 and the EFOP-3.6.2-16-2017-00015 projects. These projects are co-financed by the European Union and the European Social Fund.}

\dedicatory{Dedicated to the 75th birthday of Professor Zygfryd Kominek}

\begin{abstract}
The main goal of this paper is to give a completely elementary proof for the decomposition theorem of Wright convex functions which was discovered by C.\ T.\ Ng in 1987. In the proof, we do not use transfinite tools, i.e., variants of Rod\'e's theorem, or de Bruijn's theorem related to functions with continuous differences.
\end{abstract}

\maketitle

\section{Introduction}

In 1954, E.\ M.\ Wright \cite{Wri54} (see also \cite{RobVar73}) introduced a new convexity property for real functions: A real valued function $f$ defined on an interval $I$ is called \emph{Wright convex} if, for all $x,y\in I$ and $t\in[0,1]$,
\Eq{W}{
  f(tx+(1-t)y)+f((1-t)x+ty)\leq f(x)+f(y).
}
Here and in the sequel, $I\subseteq\R$ denotes a nonvoid open interval.

One can easily see that convex functions are Wright convex. Indeed, if $f:I\to\R$ is convex, then, for all $x,y\in I$ and $t\in[0,1]$,
\Eq{*}{
  f(tx+(1-t)y)\leq tf(x)+(1-t)f(y) \quad\mbox{and}\quad
  f((1-t)x+ty)\leq (1-t)f(x)+tf(y).
}
Adding up these two inequalities side by side, we can arrive at the conclusion that $f$ is Wright convex, too.

On the other hand, if $f:\R\to\R$ is additive, then its restriction $f|_I$ is also Wright convex on $I$. Indeed, for all $x,y\in I$ and $t\in[0,1]$, the additivity of $f$ implies
\Eq{*}{
  f(tx+(1-t)y)+f((1-t)x+ty)&=f(tx+(1-t)y+(1-t)x+ty)\\&=f(x+y)=f(x)+f(y).
}
Thus, in this case, \eq{W} holds with equality.

By the linearity of the  functional inequality \eq{W}, a function which is the sum of a convex and of an additive function must be Wright convex, too. The following surprising characterization of Wright convexity, which is basically stating the converse of the above observations, was discovered by Ng \cite{Ng87b} in 1987.

\Thm{*}{{\rm (Ng's Decomposition Theorem \cite[Corollary 5]{Ng87b})}
A function $f:I\to \R$ is Wright-convex if and only if there exist 
a convex function $g:I\to \R$ and an additive function $A:\R\to\R$ such that
\Eq{WD}{
  f(x)=g(x)+A(x) \qquad(x\in I).
}
Furthermore, the additive summand $A$ in the above decomposition can be chosen so that it vanishes on $\Q$ and, together with this condition, the decomposition \eq{WD} is unique.}

The known proofs of this important theorem are usually based on results that depend on the axiom of choice, i.e., they use transfinite induction. The original proof by Ng \cite{Ng87b} used de Bruijn's theorem \cite{Bru51} which is related to functions which have continuous differences. Another approach, which was based on Rode's theorem \cite{Rod78}, was found by Nikodem \cite{Nik89a} and also by Kominek \cite{Kom03}.

For an overview about the generalizations, stability and regularity properties of Wright convex functions, we refer to the list of references, which is possibly far from being complete but may give the impression that this field of functional equations and inequalities is still actively dealt with by many researchers.

\section{Auxiliary results}

For the new proof of Ng's decomposition theorem, we will need some auxiliary results that are contained in the following three lemmas. 

\Lem{1}{Let $D$ be a dense subset of $I$ and assume that $f:D\to\R$ is a locally Lipschitz function on $D$, i.e., for every compact subinterval $[a,b]\subseteq I$, there exists $L\geq0$ such that
\Eq{fL}{
  |f(x)-f(y)|\leq L|x-y|\qquad(x,y\in[a,b]\cap D).
}
Then $f$ admits a unique locally Lipschitz extension $g:I\to\R$.}

\begin{proof}
For a fixed element $x\in I$, let $x_n\in D$ be an arbitrary sequence converging to $x$. Then there exists a compact subinterval $[a,b]\subseteq I$ which contains each member of this sequence. Using inequality \eq{fL}, it follows that
\Eq{*}{
  |f(x_n)-f(x_m)|\leq L|x_n-x_m|\qquad(n,m\in\N).
}
The sequence $(x_n)$ being a Cauchy sequence, the above inequality implies that $(f(x_n))$ is also a Cauchy sequence, whence we get that $(f(x_n))$ is convergent.
Define $g(x)$ to be the limit of this sequence. It is easy to see that this limit does not depend on the choice of the sequence $(x_n)$, furthermore, it is also clear that $f(x)=g(x)$ if $x$ belongs to $D$ (because, then the sequence $x_n:=x$ could be taken).

To see that $g$ is locally Lipschitz on $I$, let $a,b\in I$ with $a<b$. Then there exists $L\geq0$ such that \eq{fL} holds. Let $x,y\in[a,b]$ be arbitrary and choose $x_n,y_n\in[a,b]\cap D$ such that $(x_n)$ and $(y_n)$ converge to $x$ and $y$, respectively. Then, \eq{fL}
implies that
\Eq{*}{
  |f(x_n)-f(y_n)|\leq L|x_n-y_n|\qquad(n\in\N).
}
Now, upon taking the limit as $n\to\infty$, we get
\Eq{*}{
  |g(x)-g(y)|\leq L|x-y|\qquad(x,y\in[a,b]).
}
This proves that $g$ is locally Lipschitz on $I$, indeed. The uniqueness of the function $g$ follows from the well-known fact that if two continuous real valued functions defined on $I$ coincide on a dense subset of $I$, then they coincide everywhere in $I$.
\end{proof}

\Lem{2}{
Let $f:I\cap\Q\to\R$ be a Jensen convex function, i.e., assume that $f$ satisfies the Jensen inequality on $I\cap\Q$:
\Eq{JQ}{
  f\Big(\frac{x+y}2\Big)\leq \frac{f(x)+f(y)}2\qquad (x,y\in I\cap\Q).
}
Then there exists a convex function $g:I\to\R$ such that $g|_{I\cap\Q}=f$.}

\begin{proof}
Via a standard argument (see, for instance, the proofs of Lemma 5.3.1 and Theorem 5.3.5 in \cite{Kuc85}), the Jensen convexity of $f$ on $I\cap\Q$ implies that
\Eq{CQ}{
  f(tx+(1-t)y)\leq tf(x)+(1-t)f(y)\qquad (x,y\in I\cap\Q,\, t\in [0,1]\cap\Q).
}
With the notation $u:=tx+(1-t)y$, this yields that
\Eq{JI}{
  \frac{f(x)-f(u)}{x-u}\leq \frac{f(u)-f(y)}{u-y}\qquad (x,u,y\in I\cap\Q,\,x<u<y).
}
Let $[a,b]\subseteq I$ be an arbitrary compact subinterval and fix $a',a'',b',b''\in I\cap\Q$ such that $a'<a''\leq a<b\leq b'<b''$. Then, by the repeated application of inequality \eq{JI}, for any $x,y\in [a,b]\cap\Q$ with $x<y$, we obtain
\Eq{*}{
  \alpha:=\frac{f(a')-f(a'')}{a'-a''}
   \leq\frac{f(a'')-f(x)}{a''-x}\leq\frac{f(x)-f(y)}{x-y}\leq\frac{f(y)-f(b')}{y-b'}
   \leq\frac{f(b')-f(b'')}{b'-b''}=:\beta.
}
Hence,
\Eq{*}{
  |f(x)-f(y)|\leq \max(|\alpha|,|\beta|)\cdot |x-y|\qquad (x,y\in [a,b]\cap\Q),
}
which shows that $f$ is Lipschitz with modulus $L:=\max(|\alpha|,|\beta|)$ on the set $[a,b]\cap\Q$. Therefore $f$ is locally Lipschitz on the dense set $D:=I\cap\Q$.

By applying \lem{1}, we can now conclude that $f$ admits a unique continuous extension $g:I\to\R$. Using the density of $I\cap\Q$ in $I$, the density of $[0,1]\cap\Q$ in $[0,1]$ and the continuity of $g$, the inequality \eq{CQ} easily implies that the function $g$ is also convex.
\end{proof}

\Lem{3}{Let $J\subseteq I$ be an open subinterval, let $\varphi:J\to\R$ be nondecreasing and $\psi:J\to\R$ be continuous such that $\varphi(x)=\psi(x)$ holds for all element $x$ of a dense subset $D\subseteq J$. Then $\varphi=\psi$ holds on $J$.}

\begin{proof}
Let $x\in J$ be fixed. First choose a sequence $(x_n)$ in $D$ such that $(x_n)$ converges to $x$ and $x\leq x_n$ holds for all $n\in\N$. Then, by the nondecreasingness of $\varphi$ and the continuity of $\psi$, we obtain
\Eq{*}{
  \varphi(x)
  \leq \lim_{n\to\infty}\varphi(x_n)
  = \lim_{n\to\infty}\psi(x_n)
  =\psi(x).
}
By taking another sequence $x_n\in D$ such that $(x_n)$ converges to $x$ and $x\geq x_n$ holds for all $n\in\N$, we similarly get that
\Eq{*}{
  \varphi(x)
  \geq \lim_{n\to\infty}\varphi(x_n)
  = \lim_{n\to\infty}\psi(x_n)
  =\psi(x).
}
The above two inequalities finally yield that $\varphi(x)=\psi(x)$, which was to be proved.
\end{proof}

\section{The proof of Ng's decomposition theorem}

Let $f:I\to\R$ be Wright convex. Then, with $t=\frac12$, the Wright convexity of $f$ implies that it is Jensen convex, i.e.,
\Eq{J}{
  f\Big(\frac{x+y}2\Big)\leq \frac{f(x)+f(y)}2\qquad (x,y\in I).
}
Therefore, the restriction $f|_{I\cap\Q}$ is Jensen convex on $I\cap\Q$. Thus, by \lem{2}, there exists a convex function $g:I\to\R$ such that $g|_{I\cap\Q}=f|_{I\cap\Q}$. The convexity of $g$ implies that it is also continuous on $I$.

In the rest of the proof, we show that $\Phi:=f-g$ is a restriction of an additive function to $I$ which vanishes on $\Q$. First, for all $u,v>0$, we prove that
\Eq{fg}{
  (\Delta_u\Delta_v f)(x)=(\Delta_u\Delta_v g)(x)\qquad (x\in I\cap(I-v-u)).
}
Here, for $w\in\R$, the difference operator $\Delta_w$ acts on functions $f:I\cap(I-w)\to\R$ and is defined by $(\Delta_wf)(x):=f(x+w)-f(x)$.

Observe that, with the notations $u:=t(y-x)$, $v:=(1-t)(y-x)$, the Wright convexity property can be reformulated as the inequality
\Eq{uv}{
  (\Delta_u\Delta_vf)(x)\geq0 \qquad \big(u,v>0,\, x\in I\cap(I-u-v)\big).
}
This implies that, for each fixed $v>0$, the function $\Delta_vf:I\cap(I-v)\to\R$ is nondecreasing. In particular, for a fixed $v\in\Q_+$, the function $\Delta_v f:I\cap(I-v)\to\R$ is nondecreasing, $\Delta_v g:I\cap(I-v)\to\R$ is continuous and, for $x\in I\cap(I-v)\cap\Q$, the equality of $f$ and $g$ on $I\cap\Q$ implies
\Eq{*}{
  (\Delta_v f)(x)=f(x+v)-f(x)=g(x+v)-g(x)=(\Delta_v g)(x).
}
Therefore, by applying \lem{3} to the subinterval $J:=I\cap(I-v)$, to the functions $\varphi:=\Delta_vf$, $\psi:=\Delta_vg$ and to the dense set $D:=I\cap(I-v)\cap\Q$, we obtain that the equality
\Eq{*}{
  (\Delta_v f)(x)=(\Delta_v g)(x)\qquad (x\in I\cap(I-v))
}
holds for all $v\in\Q_+$. Consequently, for all $v\in\Q_+$ and $u>0$,
\Eq{*}{
  (\Delta_u\Delta_v f)(x)=(\Delta_u\Delta_v g)(x)\qquad (x\in I\cap(I-v-u)).
}
Let $u,v\in\R_+$ be arbitrarily fixed. Then, for $v'\in\Q_+$ with $v'<v$ and for all 
$x\in I\cap(I-v-u)$, the above equality and the inequality \eq{uv} imply
\Eq{*}{
  (\Delta_v\Delta_u f)(x)&=(\Delta_u f)(x+v)-(\Delta_u f)(x)\\
   &=(\Delta_u f)(x+v)-(\Delta_u f)(x+v')+(\Delta_u f)(x+v')-(\Delta_u f)(x)\\
   &=(\Delta_{v-v'}\Delta_u f)(x+v')+(\Delta_{v'}\Delta_u f)(x)
   \geq(\Delta_{v'}\Delta_u f)(x)=(\Delta_{v'}\Delta_u g)(x).
}
Taking the limit $v'\to v$, it follows that
\Eq{*}{
  (\Delta_v\Delta_u f)(x)\geq(\Delta_{v}\Delta_u g)(x)\qquad(x\in I\cap(I-v-u)).
}
Similarly, by an analogous argument, we can also obtain
\Eq{*}{
  (\Delta_v\Delta_u f)(x)\leq(\Delta_{v}\Delta_u g)(x)\qquad(x\in I\cap(I-v-u)).
}
Therefore, for all $u,v>0$, the equality \eq{fg} holds, in other words, we have that
\Eq{*}{
  (\Delta_v\Delta_u (f-g))(x)=0\qquad(x\in I\cap(I-v-u)).
}
Thus, for $x,y\in I$ with $x<y$, the above equality with $u:=v:=\frac{y-x}{2}$ implies that $\Phi:=f-g$ satisfies the Jensen equation, i.e., 
\Eq{*}{
    \Phi\Big(\frac{x+y}2\Big)=\frac{\Phi(x)+\Phi(y)}2.
}
By classical results on this equation (cf.\ \cite{Kuc85}), $\Phi$ is of the form $\Phi(x)=A(x)+c$\, $(x\in I)$, where $A:\R\to\R$ is an additive function and $c\in\R$ is a constant. Since $\Phi|_{I\cap\Q}=(f-g)|_{I\cap\Q}=0$, we get that $A(1)x+c=A(x)+c=\Phi(x)=0$ for all $x\in I\cap\Q$. This implies that $A(1)=c=0$ and hence $\Phi$ equals to the restriction of an additive function which vanishes on $\Q$. Thus the proof is completed.

\def\MR#1{}

\end{document}